\documentstyle[art10]{article}
\title{Classification of Elliptic Line Scrolls}
\author{Luis Fuentes\thanks{Supported by an F.P.I.
fellowship of Spanish Government}}
\date{}

\newtheorem{teo}{Theorem}[section]

\newtheorem{prop}[teo]{Proposition}

\newtheorem{lemma}[teo]{Lemma}

\newtheorem{rem}[teo]{Remark}

\def\vhi{\varphi}

\def\p{{\bf P}}
\def\E{{\cal E}}
\def\Te{{\cal O}}

\font\euf=eufm10 at 12pt
\def\e{\mbox{\euf e}}
\def\b{\mbox{\euf b}}
\def\aa{\mbox{\euf a}}
\def\A{{\alpha}}

\def\P{{\bf P}}
\def\gr{\partial}

\def\qed{\hspace{\fill}$\rule{2mm}{2mm}$}

\def\lrw{{\longrightarrow}}
\def\rw{{\rightarrow}}

\begin{document}

\maketitle

{\footnotesize{\bf Author's address:} Departamento de Algebra, Universidad de Santiago
de Compostela. $15706$ Santiago de Compostela. Galicia. Spain. e-mail: {\tt fuentes@zmat.usc.es}\\
{\bf Abstract:} The classification of projective elliptic line scrolls with the
description of their singular loci is given. In particular we recover Atiyah
Theorem by using classical methods.}

\vspace{0.1cm}

{\bf Introduction:} Through this paper, a {\it geometrically ruled surface}, or simply a {\it
ruled surface}, will be a $\P^1$-bundle over a smooth curve $X$ of genus $g$. It will be
denoted by $\pi: S=\P(\E_0)\lrw X$ and we will follow the notation and
terminology of R. Hartshorne's book \cite{hartshorne}, V, Section 2. We will suppose that
$\E_0$ is a normalized sheaf and $X_0$ is the section of minimum self-intersection
that corresponds to the surjection $\E_0\lrw \Te_X(\e)\lrw 0$, $\bigwedge^2\E\cong
\Te_X(\e)$. Consider the following question: which are the linear equivalence classes $D\sim mX_0+\b
f$,
$\b\in Pic(X)$, that correspond to very ample divisors? When $g=1$ and $m=1$, a characterization is
known (\cite{hartshorne}, V, Ex. 2.12), but the classification of elliptic scrolls obtained
by Corrado Segre in \cite{corrado} does not follow directly from this. A scroll is the
image of a ruled surface $\pi:S=\P(\E_0)\lrw X$ by a unisecant complete linear
system.

In this paper, I apply the results of the previous work \cite{fuentes} to obtain the classification of
elliptic (line) scrolls in $\P^N$, $N\geq 3$, including the singular ones and degenerations. Theorems
\ref{scrollselipticasdescomp} and \ref{scrollselipticasnodescomp} provide the description of all
scrolls in the decomposable and indecomposable cases respectively. 

To obtain Theorem \ref{scrollselipticasnodescomp}, we use that all indecomposable ruled surface is
obtained from a decomposable one by applying a finite number of elementary transformations, see
(\cite{fuentes}, 3.9). We give a proof of Atiyah Theorem (the unique indecomposable elliptic scrolls
have $e=0$ or $e=-1$) that seems geometrically interesting to us.

In Theorem \ref{telementalelipticasdescomp}, we prove that the elementary transformation at a point of
an elliptic decomposable ruled surface is either decomposable or indecomposable with $e=0,-1$.
Proposition
\ref{telementalnodescomp0} says that an indecomposable elliptic ruled surface with $e=0$ is obtained
from a decomposable one or from an indecomposable one with $e=-1$. Finally, in Proposition
\ref{telementalnodescomp1} we see that the elementary transform at a point of an indecomposable
elliptic ruled surface with $e=-1$ is either a decomposable one or an indecomposable one with $e=0$.
This result uses Proposition \ref{irreduciblesnodescomp1}, where we give a parameterization of the
indecomposable ruled surface with $e=-1$ that makes it isomorphic to the symmetric product
$S^2X$. Moreover, this result provides explicitly how an elliptic ruled surface is obtained from
$X\times P^1$ by applying elementary transformations. Compare with \cite{maruyama1}, II, $\S2$.

We give tables that show the classification of elliptic line scrolls in $\P^N$, $N\geq 3$, and contain
the results of C. Segre in \cite{corrado}.

Finally, in $\S3$ we study the classes of $2$-secant divisors in a elliptic ruled surface. We
characterize those that are base-point-free and very ample, and those that contain irreducible
elements. These results make exhaustive our study of elliptic scrolls.

In a forthcoming paper we will study the classification of elliptic scrolls of higher rank.

I would like to thank Manuel Pedreira for all his help and constant encouragement.
  
\section{Decomposable elliptic ruled surfaces.}\label{regladaselipticasdescomp}

We are going to study elliptic geometrically ruled surfaces; that is,
ruled surfaces over a smooth curve $X$ of genus $1$.

We begin by working with decomposable elliptic ruled surfaces, and we apply the results seen
in (\cite{fuentes}, section 3). Then, by using elementary transformations, we will study the
indecomposable ones.

Let us first review some properties of divisors on an elliptic curve. We will apply
them to study elliptic ruled surfaces.

\begin{lemma}\label{curvaeliptica0}

Let $X$ be an elliptic curve. Given a point $P_0\in X$, there exists an one-to-one
correspondence between points $P$ of $X$ and divisors $\b$ of $Pic_0(X)$. Under this
correspondence $\b\sim P-P_0$.\qed

\end{lemma}

\begin{lemma}\label{curvaeliptica}

Let $X$ be a smooth elliptic curve and let $\b$ be a divisor on $X$. Then:
\begin{enumerate}

\item $\b$ is base-point-free if and only if $\b\sim 0$ or $deg(\b)\geq 2$.

\item $\b$ is very ample if and only if $deg(\b)\geq 3$.\qed

\end{enumerate}

\end{lemma}

\begin{prop}\label{pfijoselipticasdescomp}

Let $S$ be a decomposable elliptic ruled surface and let $|H|=|X_0+\b f|$ be a
complete unisecant linear system. Then:

\begin{enumerate}

\item If $\e\sim 0$, then $|H|$ is base-point-free if and only if $deg(\b)\geq 2$ or
$\b\sim 0$.

\item If $\e\not\sim 0$, then $|H|$ is base-point-free if and only if $deg(\b)\geq
e+2$ or $\b\sim -\e$ and $e\geq 2$.

\end{enumerate}

\end{prop}
{\bf Proof:}

\begin{enumerate}

\item Let us suppose $\e\sim 0$. By Proposition $2.3$ in \cite{fuentes}, $|H|$ is
base-point-free if and only if $\b$ is base-point-free, because $\b+\e\sim \b$ in this
case. Applying Lemma \ref{curvaeliptica} the conclusion follows.

\item Let us suppose $\e\not\sim 0$. $|H|$ is base-point-free if and only if $\b$ and
$\b+\e$ are base-point-free. By applying Lemma \ref{curvaeliptica} we see that there
are the following possibilities: 

\begin{enumerate}

\item $\b\sim 0$ and $\b+\e\sim 0$, but then $\e\sim 0$, which contradicts our
assumption.

\item $\b\sim 0$ and $deg(\b+\e)\geq 2$, but then $deg(\e)\geq	2$ and by
(\cite{hartshorne}, V, 2.12), $S$ is indecomposable, which is false by hypothesis.

\item $deg(\b)\geq 2$ and $\b+\e\sim 0$, then $\b\sim -\e$ and necessarily
$e=deg(-\e)\geq 2$.

\item $deg(\b)\geq 2$ and $deg(\b+\e)\geq 2$, but since $deg(\b)\leq 0$, it is 
enough that $deg(\b)\geq 2+e$. \qed

\end{enumerate}

\end{enumerate}

\begin{prop}\label{amplitudelipticasdescomp}

Let $S$ be a decomposable elliptic ruled surface and let $|H|=|X_0+\b f|$ be a
complete unisecant linear system on $S$. Then $|H|$ is very ample if and only if
$deg(\b)\geq 3+e$.

\end{prop}
{\bf Proof:}

By Theorem $2.8$ in \cite{fuentes}, $|H|$ is very ample if and only if $\b$ and
$\b+\e$ are very ample. Applying Lemma \ref{curvaeliptica}, we see that this condition
holds when $deg(\b)\geq 3$ and $deg(\b+\e)\geq 3$. Since $e\geq 0$, it is sufficient
that $deg(\b)\geq 3+e$. \qed

\begin{prop}\label{irreducibleselipticasdescomp}

Let $S$ be a decomposable elliptic ruled surface. The unisecant linear systems with
generic element irreducible are:

\begin{enumerate}

\item If $\e\sim 0$, they are $|X_0|$ and $|X_0+\b f|$ with $deg(\b)\geq 2$.

\item If $\e\not\sim 0$, they are $|X_0|$, $|X_0-\e f|$ and $|X_0+\b f|$ with
$deg(\b)\geq 1+e$.

\end{enumerate}

\end{prop}
{\bf Proof:}

We will apply Theorem $2.5$ in \cite{fuentes}. The linear system $|X_0+\b
f|$ have irreducible elements if and only if either $\b\sim 0$, or $\b\sim -\e$, or $\b$ and $\b+\e$
are effective without common base points. Thus we have two cases:

\begin{enumerate}

\item If $\e\sim 0$, then $\b+\e\sim \b$. So $|X_0+\b f|$ has irreducible elements if
$\b\sim 0$ or $\b$ is base-point-free. According to Lemma \ref{curvaeliptica}, we know
that this happens when $\b\sim 0$ or $\deg(\b)\geq 2$.

\item If $\e\not\sim 0$, we have irreducible elements in $|X_0|$, $|X_0-\e f|$ and $|
X_0+\b f|$ when $\b$ and $\b+\e$ are effective without common base points.

$\b+\e$ is effective if $deg(\b+\e)\geq 0$, this is, if $deg(\b)\geq e$. Moreover, if
$deg(\b)=e$, then $\b\sim -\e$ and the system is $|X_0-\e f|$. Thus, if $deg(\b)\geq 1+e$,
then $\b$ and $\b+\e$ are effective and they can only have base points when $e=0$. In
this case, since $e\not\sim 0$ and $deg(\b)=1$, $\b$ and $\b+\e$ correspond to
different points of $X$. So they have not common base points. \qed

\end{enumerate}

\begin{lemma}\label{normalidadelipticasdescomp}

Let $S$ be a decomposable elliptic ruled surface, let $|H|=|X_0+\b f|$ be a complete
base-point-free linear system on $S$ such that $\b\not\sim 0$, and let $\phi:S\lrw \P^N$ be the
regular map defined by $|H|$. Let $D$ be an irreducible
unisecant curve on $S$ that is not linearly equivalent to $X_0$ and $X_1$, with $D\sim X_0+\aa f$.
Then:

\begin{enumerate}

\item If $\b\not\sim -\e$, $\phi(D)$ is linearly normal if and only if $deg(\aa)\leq
deg(\b)$ and $\aa\not\sim \b$.

\item If $\b\sim -\e$, $\phi(D)$ is linearly normal if and only if $deg(\aa)=1+e$.

\end{enumerate}

\end{lemma}
{\bf Proof:}

Let us consider the trace of the complete linear system $|H|$ on the curve $D$:
$$
\begin{array}{l}
0\rw H^0(\Te_S((\b -\aa) f))\lrw
H^0(\Te_S(X_0+\b f))\stackrel{\A}{\lrw} H^0(\Te_{D}(X_0+\b
f))=\\
=H^0(\Te_X(\aa+\b+\e))\rw 0
\end{array}
$$
The curve $\phi(D)$ is linearly normal when $|H|$ traces the complete linear system
$|\aa+\b+\e|$ on $D$; equivalently, when $\A$ is a surjection ($h^0(\Te_S(X_0+\b
f))=h^0(\Te_X(\aa+\b+\e))+h^0(\Te_X(\b-\aa))$). 

We know that $h^0(\Te_S(X_0+\b f))=h^0(\Te_X(\b))+h^0(\Te_X(\b+\e))$.

\begin{enumerate}

\item Let us first suppose that $\b\not\sim -\e$.

$D$ is not linearly equivalent to $X_0$ or $X_1$. So by 
Proposition \ref{irreducibleselipticasdescomp}, $deg(\aa+\e)\geq 1$. Moreover,  $|H|$
is base-point-free and, by hypothesis, $\b\not\sim 0$ and $\b\not\sim -\e$; by applying
Proposition~\ref{pfijoselipticasdescomp} we have that 
$deg(\b)\geq 2+e$ and then
$deg(\aa+\e+\b)\geq 1$ and the divisor $\aa+\e+\b$ is nonspecial. By Riemann-Roch,
$h^0(\Te_X(\aa+\b+\e))=deg(\aa)+deg(\b)+deg(\e)$ and
$h^0(\Te_X(\b-\aa))=deg(\b)-deg(\aa)+h^1(\Te_X(\b-\aa))$.

Since $deg(\b)\geq 2+e$, $\b$ and $\b+\e$ are nonspecial on
$X$ and by Riemann-Roch, $h^0(\Te_S(X_0+\b
f))=2deg(\b)-e$. Therefore, $h^0(\Te_S(X_0+\b
f))-(h^0(\Te_X(\aa+\b+\e))-h^0(\Te_X(\b-\aa)))=-h^1(\Te_X(\b-\aa))$,
and the map $\A$ is a surjection when $\b-\aa$ is nonspecial. That happens
when $\b\not\sim \aa$ and $deg(\aa)\leq deg(\b)$.

\item Let us suppose that $\b\sim -\e$. 

Since $|X_0+\b f|$ is base-point-free, $h^0(\Te_S(X_0+\b
f))=h^0(\Te_X(\b))+h^0(\Te_X(\b+\e))=e+1$.

Moreover, $D$ is not linearly equivalent to $X_0$ or $X_1$. So by
Proposition~\ref{irreducibleselipticasdescomp}, we have that $deg(\aa+\e)\geq 1$ (equivalently,
$deg(\aa)\geq 1+e$).

Thus, 
$h^0(\Te_X(\aa+\b+\e))=h^0(\Te_X(\aa))=deg(\aa)$ and
$h^0(\Te_X(\b-\aa))=h^0(\Te_X(\e-\aa))=0$; from this, the equality
$$h^0(\Te_S(X_0+\b f))=h^0(\Te_X(\aa+\b+\e)+h^0(\Te_X(\b-\aa))$$
holds if and only if $deg(\aa)=1+e$.
\qed

\end{enumerate}

\bigskip

Finally, a direct application of the above results allows us to describe decomposable
elliptic scrolls in $\P^N$.

\begin{teo}\label{scrollselipticasdescomp}

Let $S$ be a decomposable elliptic ruled surface. Let $|H|=|X_0+\b f|$ be a complete
base-point-free linear system on $S$, with $b=deg(\b)$. Let $\phi:S\lrw \P^N$ be the
regular map defined by $|H|$ on $S$ and let $R=\phi(S)$ be the image scroll. Then
$R$ is one of the following models:

\begin{enumerate}

\item When $\phi$ is not birational:

\begin {enumerate}

\item If $\e\sim 0$ and $\b\sim 0$, then $R$ is a line parameterizing the curves of
$|X_0|$. In fact, $S\cong X\times \P^1$ and $\phi$ is the second projection.

\item If $\e\sim 0$ and $b=2$, then $R$ is a smooth quadric in $\P^3$ and $\phi$ is a
2:1 morphism from $S$ onto $R$.

\item If $e=2$ and $\b\sim -\e$, then $R$ is a plane and $\phi$ is a 2:1 morphism.

\end{enumerate}

\item When $\phi$ is birational:

\begin{enumerate}

\item If $e>2$ and $\b\sim -\e$, then $R$ is a cone in $\p^e$ over a smooth linearly
normal elliptic curve of degree $e$.
 
The singular locus of $R$ is the vertex of the cone.
 
$R$ has a
family of linearly normal elliptic curves of degree $e$ that correspond to the hyperplane sections.
 
$R$
has unisecant elliptic curves of degree $d>e$ in the linear systems $|X_0+\aa f|$ with $deg(\aa)=d$.
Theses curves are smooth and linearly normal if and only if
$d=1+e$. If $d>1+e$ then the curves have a singular point of multiplicity $d-e$ at
the vertex of the cone.

\item If $e=0$, $\e\not\sim 0$ and $b=2$, then $R$ is a elliptic scroll of degree $4$
in $\P^3$. It is generated by a $2:2$ correspondence between two disjoint lines. 

The singular locus of
$R$ are disjoint lines which generate it. 

$R$ has unisecant elliptic curves of degree $d\geq 3$ in the
linear systems $|X_0+\aa f|$ with $deg(\aa)=d-2$. These curves are linearly normal if and only if
$d\leq 4$ and 
$\aa\not\sim \b$ (if $\aa\sim \b$ they are the hyperplane sections).

\item If $e>0$ and $b=e+2$, then $R$ is a elliptic scroll of degree $e+4$ in
$\P^{e+3}$. It is generated by a $1:2$ correspondence between a line and a smooth
linearly normal elliptic curve of degree $e+2$, laying in disjoint spaces. 

The singular locus of $R$ is
the directrix line. 

$R$ has a unique directrix curve of minimum degree $2$ and unisecant elliptic
curves of degree $d\geq 3+e$ in the linear systems $|X_0+\aa f|$ with $deg(\aa)=d-2$. These
curves are linearly normal if and only if $d\leq 4+e$ and 
$\aa\not\sim \b$ (if $\aa\sim \b$ they are the hyperplane sections).

\item If $b\geq 3+e$, then $R$ is a smooth elliptic scroll of degree $2b-e$ in
$\p^{2b-e-1}$. It is generated by a $1:1$ correspondence between two smooth linearly
normal elliptic curves of degrees $b$ and $b+e$, laying in disjoint spaces. 

If $\e\sim 0$, $R$ has a
one-dimensional family of disjoint directrix curves of minimum degree $b$. Moreover, $R$ has unisecant
elliptic curves of degree $d\geq b+2$ in the linear systems $|X_0+\aa f|$ with $deg(\aa)=d-b$. These
curves are linearly normal if and only if $d\leq 2b$ and 
$\aa\not\sim \b$ (if $\aa\sim \b$ they are the hyperplane sections). 

If $\e\not\sim 0$, $R$ has an
directrix curve of minimum degree $b-e$ and unisecant linearly normal elliptic curves of degree $b$ in
the linear system $|X_0-\e f|$. Therefore, if $e>0$ the minimum degree curve is unique, but if $e=0$
and $\e\not\sim 0$, there are two curves of minimum degree. Moreover,
$R$ has unisecant elliptic curves of degree $d\geq b+1$ in the linear systems
$|X_0+\aa f|$ with $deg(\aa)=d+e-b$. These curves are linearly normal if and only if
$d\leq 2b-e$ and 
$\aa\not\sim \b$ (if $\aa\sim \b$ they are the hyperplane sections).

\end{enumerate}

\end{enumerate}

\end{teo}
{\bf Proof:}

Note that, when $\phi$ is birational, the description of the scroll $R$ follows
immediately from Theorem $2.9$ in \cite{fuentes}. The families of irreducible
unisecant curves are described in Proposition \ref{irreducibleselipticasdescomp} and we know
when they are linearly normal by Lemma \ref{normalidadelipticasdescomp}. If
$|H|$ is very ample, then $\phi$ is an isomorphism. Let us study other cases.

By Propositions \ref{pfijoselipticasdescomp} and \ref{amplitudelipticasdescomp}, we
know when $|H|$ is base-point-free but not very ample:

\begin{enumerate}

\item If $e\sim 0$ and $\b\sim 0$, then $|H|=|X_0|$. The linear system has dimension
one, so $R$ is a line parameterizing curves of $|X_0|$. Each curve of this linear
system is isomorphic to $X$ and it is applied onto a point of $R$. From this, the
surface
$S$ is isomorphic to $X\times \P^1$ and $\phi$ is the second projection.

\item If $\e\sim 0$ and $b=2$, then, by Theorem $2.9$ in \cite{fuentes}, $\phi(X_0)$ and
$\phi(X_1)$ are two disjoint lines given by the $2:1$ morphism defined by $|\b|$ from
$X$ onto $\P^1$. Moreover, $P$ and $Q$ apply on the same point of $\P^1$ when $\b\sim
P+Q$, this is, when $Q$ is a common base point of $\b-P$ and $\b+\e-P$. Hence,
generators $Qf$ and $Pf$ apply onto the same line on $R$ if $P+Q\sim \b$.
Consequently, a unique (double) line passes through each point of $\phi(X_0)$ and
$\phi(X_1)$. $R$ is generated by a $1:1$ correspondence between two disjoint lines in
$\P^3$, so $R$ is the smooth quadric of $\P^3$ and $\phi$ is a $2:1$ morphism.

\item If $e=2$ and $\b\sim -\e$, then $R$ have degree $2$ in $\P^3$. Since $\b+\e\sim
0$ and $b=2$, $\phi(X_0)$ is a point and $\phi(X_1)$ is a double line. The scroll is
generated by lines meeting the point and the double line. Therefore, $R$ is a plane and
$\phi$ is a $2:1$ map.

\item If $e>2$ and $\b\sim -\e$, since $\b+\e\sim 0$, $\phi(X_0)$ is a point which meets
all generators. Then $R$ is a cone with vertex $\phi(X_0)$ over the linearly normal
elliptic curve $\phi(X_1)$. 

By Proposition \ref{irreducibleselipticasdescomp} we know that $S$ have unisecant
irreducible curves in the linear systems $|X_0|$ (the vertex of the cone) , $|X_1|$
(the hyperplane sections) and $|X_0+\aa f|$ with $deg(\aa)\geq 1+e$. These curves
have degree $d=(X_0+\aa f).(X_0-\e f)=d$. If $deg(\aa)\geq 1+e$, the curves meets
$X_0$ at $X_0.(X_0+\aa f)=deg(\aa)-e\geq 2$ points. So they have a singular point of
multiplicity $d-e$ at the vertex $\phi(X_0)$.
 
\item If $\e\not\sim 0$ and $b=e+2$, then $\b-P$ and $\b+\e-P$ have
not common base points for any $P\in X$. So $\phi$ is birational. By applying Theorem $2.9$ in
\cite{fuentes}, our claims follow. \qed

\end{enumerate}

\section{Indecomposable elliptic ruled surfaces.}\label{regladaselipticasnodescomp}

Our goal is to find models of indecomposable elliptic ruled surfaces. According to
Corollary $3.9$ in \cite{fuentes}, they are obtained from decomposable ones by
applying a finite number of elementary transformations.

In the following theorem we apply Theorem $3.10$ in \cite{fuentes}, to study how
a decomposable ruled surfaced is modified by an elementary transformation.

\begin{teo}\label{telementalelipticasdescomp}

Let $\pi:S\lrw X$ be a decomposable elliptic ruled surface. Let $x\in Pf$ be a point
of $S$. Let $S$ be the elementary transform of $S'$ at $x$ corresponding to a
normalized sheaf $\E_0'$ with divisor $\e'$, $e'=-deg(\e')$. Let $Y_0$ be the
minimum self-intersection curve of $S'$. Then we have the following cases:

\begin{enumerate}

\item When $e\geq 2$:

\begin{enumerate}

\item If $x\in X_0$, then $S'$ is decomposable, $\e'\sim \e-P$ and $Y_0=X_0'$.

\item If $x\notin X_0$, then $S'$ is decomposable, $\e'\sim \e+P$ and $Y_0=X_0'$.

\end{enumerate}

\item When $e=1$:

\begin{enumerate}

\item If $x\in X_0$, then $S'$ is decomposable, $\e'\sim \e-P$ and $Y_0=X_0'$.

\item If $x\notin X_0$ and $-\e\not\sim P$, then $S'$ is decomposable, $\e'\sim \e+P$ and
$Y_0=X_0'$.

\item If $x\notin X_0$, $x\in X_1$ and $\e\sim -P$, then $S'$ is decomposable, $\e'\sim 0$ and
$Y_0=X_0'$.

\item If $x\notin X_0$, $x\notin X_1$ and $\e\sim -P$, then $S'$ is indecomposable, $\e'\sim 0$
and $Y_0=X_0'$.

\end{enumerate}

\item When $e=0$ and $\e\not\sim 0$:

\begin{enumerate}

\item If $x\in X_0$, then $S'$ is decomposable, $\e'\sim \e-P$ and $Y_0=X_0'$.

\item If $x\in X_1$, then $S'$ is decomposable, $\e'\sim -\e-P$ and $Y_0=X_1'$.

\item If $x\notin X_0$, then $x\notin X_1$, $S'$ is indecomposable, $\e'\sim \e+P$ and
$Y_0=X_0'$.

\end{enumerate}

\item If $\e\sim 0$, then $S'$ is decomposable, $\e'\sim -P$ and $Y_0=X_0'$.

\end{enumerate}

\end{teo}
{\bf Proof:}

It is sufficient to apply Theorem $3.10$ in \cite{fuentes},  directly:

\begin{enumerate}

\item When $e\geq 2$:

If $x\in X_0$ we apply the point one of Theorem $3.10$ in \cite{fuentes}.

If $x\notin X_0$, since $e\geq 2$, $h^0(\Te_X(-\e))>0$ and $-\e$ has not base points.
By point $3$ of Theorem $3.10$ in \cite{fuentes}, $S'$ is decomposable.

\item If $e=1$, then $h^0(\Te_X(-\e))=1>0$ and $-\e$ has a unique base point.
Applying the four points of Theorem $3.10$ in \cite{fuentes}, we deduce four assertions.

\item If $e=0$ and $\e\not\sim 0$, then $h^0(\Te_X(-\e))=0$ and any point of $X$ is a
base point of $-\e$. Using points 1,2 and 3 of Theorem $3.10$ in \cite{fuentes}, the
conclusion follows.

\item If $\e\sim 0$, then $h^0(\Te_X(-\e))=1$ and $-\e$ has not base points. By
points 1 and 3 of Theorem $3.10$ in \cite{fuentes}, the 
ruled surface $S'$ is always decomposable, with $\e'\sim -P$ and $Y_0=X_0'$. \qed

\end{enumerate}

\bigskip

We have obtained two models of indecomposable elliptic ruled surfaces by applying an
elementary transformation to a decomposable one.

If we repeat this construction with these models, there could appear new
indecomposable surfaces. However, we will see that there are no other models of
indecomposable elliptic ruled surfaces.

In order to get how these ruled surfaces are modified by elementary transformations,
we begin by investigating their families of unisecant curves.

\begin {prop}\label{irreduciblesnodescomp0}

Let $\pi:S\lrw X$ be an indecomposable elliptic ruled surface with invariant $e=0$ and
$\e\sim 0$. The complete linear system $|H|=|X_0+\b f|$ has irreducible elements if and
only if $\b\sim 0$ or $deg(\b)\geq 1$. Moreover, $X_0$ is the unique minimum
self-intersection curve.

\end{prop}
{\bf Proof:}

Let us see that $X_0$ is the unique minimum self-intersection curve. Because $e=0$,
if there existed a curve $C$ with $C^2=0$, then $C$ would not meet $X_0$. So the
surface would have two disjoint unisecant curves and it would be decomposable.

Let $D\sim X_0+\b f$ be an irreducible element different from $X_0$. It satisfies
$\pi_*(X_0\cap D)\sim \e+\b\sim \b$. Then, $\b$ must be effective and $deg(\b)\geq 1$.

Conversely, suppose $deg(\b)\geq 1$:

\begin{enumerate}

\item If $deg(\b)\geq 2$, then $\b$ is nonspecial. Moreover, $\b$ and $\b+\e$ have no
common base points. So, by Corollary $1.8$ in \cite{fuentes}, the linear system has
irreducible elements.

\item If $deg(\b)=1$, then, 
$$h^0(\Te_S(X_0+\b f))=h^0(\Te_X(\b))+h^0(\Te_X(\b+\e))=2$$
because $\b$ is nonspecial. Reducible elements of $|X_0+\b f|$ will contain generators
and they will left a residual unisecant curve. As $\deg(\b)=1$, this curve must have
self-intersection $0$, and so it is $X_0$. From this, reducible elements of $|H|$ are in the
linear subsystem $|H-X_0|=|\b f|$. Since $h^0(\Te_S(\b f))=1<h^0(\Te_S(X_0+\b f))$, the
generic element of $|H|$ is irreducible. \qed

\end{enumerate}

\begin {prop}\label{irreduciblesnodescomp1}

Let $\pi:S\lrw X$ be an indecomposable elliptic ruled surface with invariant $e=-1$ and
$\e\sim P$. The complete linear system $|H|=|X_0+\b f|$ has irreducible elements if and
only if $\deg(\b)\geq 0$. Moreover, there is an unidimensional family of curves of minimum
self-intersection parameterized by the base elliptic curve $X$. Two of these curves pass through any
point of each generator $Qf$, except through four points $\{x_i,1\leq i\leq 4\}$. Through
these points it passes a unique a curve $D_i$, with $D_i\sim X_0+(R_i-P)f$. The points
$R_i$ are the ramification points of the morphism $|Q+P|:X\lrw \P^1$; that is, points
satisfying $2R_i\sim Q+P$.

\end{prop}
{\bf Proof:}
Any unisceant curve $D\sim X_0+\b f$ has self-intersection greater than or equal to
$X_0^2$, so necessarily $deg(\b)\geq 0$.

If $deg(\b)\geq 1$, then $\b$ and $\b+\e$ are effective divisors, $\b$ is nonspecial
and $\b+\e$ is base-point-free. By applying Corollary $1.8$ in \cite{fuentes}, we see
that the generic element of $|X_0+\b f|$ is irreducible.

In Lemma \ref{curvaeliptica0} we saw that, given a point $P\in X$, a divisor $\b$ of
degree $0$ can be written as $\b\sim Q-P$. In this way, points $Q$ of $X$ parameterize
divisors of degree $0$ on $X$.

Let $\b\sim Q-P$ be a divisor of degree $0$ with $Q\neq P$. Then $h^0(\b)=0$ and
$\b$ is nonspecial. According to Remark $1.2$ in \cite{fuentes}, we see that
$$h^0(\Te_S(X_0+\b f))=h^0(\Te_X(\b))+h^0(\Te_X(\b+\e))=1.$$
Because $deg(\b-R)<0$, then $h^0(\Te_S(X_0+(\b-R) f))=0$
for any point $R\in P$. Therefore, by Proposition $1.7$ in
\cite{fuentes},
$|X_0+(Q-P)|$ contains a unique irreducible curve that will be denoted by $D_Q$.

Summarizing, we have a family of curves $D_Q$ with self-intersection $1$ par\-ameterized
by $X$. If $Q\neq P$, then we know that $dim(|D_Q|)=0$. Let us see that there is a
unique curve in the linear
system $|X_0|=|D_P|$. Considering the trace of $|X_0|$ on
$|D_Q|$ with $P\neq Q$, we find:
$$
0\rw H^0(\Te_S((P-Q) f))\rw H^0(\Te_S(X_0))\stackrel{\A}{\lrw}
H^0(\Te_{D_Q}(Q))\rw H^1(\Te_S((P-Q)f))
$$

Since $P-Q\not\sim 0$, then $h^1(\Te_S((P-Q)f))=0$ and
$$h^0(\Te_S(X_0))=h^0(\Te_X(P-Q))+h^0(\Te_X(Q))=1$$.

We now study how curves of this family intersect. It is clear that $D_Q.D_R=1$.
Moreover, $\pi_*(D_Q\cap D_R)\sim Q+R-P$; consequently, two curves $D_Q$ and $D_R$
meets at a point on the generator $Tf$, with $T$ satisfying $T+P\sim Q+R$.

In this way, given a generator $Tf$ we define the map $\sigma:X\lrw Tf\cong P^1$ by
assigning $D_Q\cap Tf$ to each point $Q$ of $X$.

Let us fix a curve $D_Q$. Other curve $D_R$ meets $D_Q$ on $Tf$ when $Q+R\sim T+P$;
that is, when $R\sim T+P-Q$. Hence, there is a unique curve meeting $D_Q$ on $Tf$.
This curve coincides with $D_Q$ when $2Q\sim T+P$.

We have seen that the morphism $\sigma$ is not constant and it applies two points $Q$
and $R$ satisfying $Q+R\sim T+P$ onto a point of $\P^1$. The morphism $\sigma$ is defined by the
linear system $|T+P|$ on $X$. It is a 2:1 morphism from a elliptic curve onto a line.
By Hurwitz Theorem (\cite{hartshorne}, IV,2.4), $\sigma$ has exactly four
ramifications $\{R_i,1\leq i\leq 4\}$ satisfying $2R_i\sim P+Q$. It follows
that a unique curve of the family passes through them. \qed

\bigskip

\begin{rem}\label{parametrizacion}
The above proposition allows to obtain a parameterization of the indecomposable elliptic ruled surface
$\P(\E_0)$ with $e=-1$ and $\e\sim P$ in the following way.

Let $S^2X$ be the divisors of degree $2$ in $X$. We define the map:
$$
\begin{array}{rcl}
{S^2X}&{\stackrel{\tau}{\lrw}}&{\P(\E_0)}\\
{Q+R}&{\lrw}&{D_Q\cap D_R}\\
\end{array}
$$
The map is well defined except at most at points of the diagonal of $S^2X$. But we have seen that the
curves $D_Q, D_R$ meet at a point in the generator $Tf$ such that $Q+R\sim T+P$. Thus, equivalently,
the map can be defined:
$$
\begin{array}{rcl}
{S^2X}&{\stackrel{\tau}{\lrw}}&{\P(\E_0)}\\
{Q+R}&{\lrw}&{D_Q\cap Tf/Q+R\sim P+T}\\
\end{array}
$$
Therefore, the map is well defined in whole $S^2X$ and it is clearly an isomorphism.

The image of the diagonal of $S^2X$ by the map $\tau$ is the focal curve $C_f$ in $\P(\E_0)$. It
meets each generator $Tf$ precisely at the four points corresponding to the ramifications of the
morphism
$\vhi_{|T+P|}:X\lrw \P^1$. By the above proposition we know that, through these points, it passes a
unique curve of the family of minimum self-intersection curves.
\end{rem}

Let us now see what happens when we apply an elementary transformation to models of
indecomposable elliptic ruled surfaces. 

\begin{prop}\label{telementalnodescomp0}

Let $\pi:S\lrw X$ be an indecomposable elliptic ruled surface with invariant $e=0$ and
$\e\sim 0$. Let $x$ be a point of $X$ with $\pi(x)=P$. Let $S'$ be the elementary
transformation of $S$ at $x$ corresponding to a normalized sheaf $\E_0'$ with
invariant $\e'$, $e'=-deg(\e')$. Let $Y_0$ be the minimum self-intersection curve of
$S'$. Then:

\begin{enumerate}

\item If $x\in X_0$, then $S'$ is decomposable with $\e'\sim -P$ and $Y_0=X_0'$.

\item If $x\notin X_0$, then $S'$ is indecomposable with $\e'\sim P$ and $Y_0=X_0'$.

\end{enumerate}

\end{prop}
{\bf Proof:}

\begin{enumerate}
\item If $x\in X_0$, then, by applying Theorem $3.8$ in \cite{fuentes}, we obtain
that $\e'\sim \e-P\sim -P$ and $Y_0=X_0'$. Moreover, by Proposition
\ref{irreduciblesnodescomp0}, we know that there exists an irreducible curve $D\sim
X_0+\P f$. But $\pi_*(D\cap X_0)\sim P$, so $D$ meets $X_0$ at $X_0\cap Pf=x$.
Since $x\in X_0$ and $x\in D$, according to elementary transformation properties
we know that $D'X_0'=DX_0-1=0$. From this, $S'$ has two disjoint unisecant curves so
it is decomposable.

\item Let us suppose $x\notin X_0$. We know $X_0$ is the unique curve of
self-intersection $0$. Any other unisecant curve $D$ of $S$ is in a linear system
$|X_0+\b f|$ with $deg(\b)\geq 1$, so $D^2\geq 2$.
By applying elementary transformation properties, we conclude $X_0'^2=X_0^2+1=1$. For any
other curve
$D$ we know that $D'^2\geq D^2-1\geq 1$. Therefore, $X_0'$ is the minimum
self-intersection curve of $S'$ and $\e'\sim \e+P$.
Finally, since $deg(\e')>0$, the ruled surface $S'$ is indecomposable. \qed

\end{enumerate}

\begin{prop}\label{telementalnodescomp1}

Let $\pi:S\lrw X$ be an indecomposable elliptic ruled surface with invariant $e=-1$ and
$\e\sim P$. Let $x$ a point of $S$ with $\pi(x)=Q$. Let $S'$ be the elementary transform
of $S$ at $x$ corresponding to a normalized sheaf $\E_0'$ with invariant $\e'$, with
$e'=-deg(\e')$ and let $Y_0$ be the minimum self-intersection curve of $S'$. Then:

\begin{enumerate}

\item If $x\in D_R$ with $D_R\sim X_0+(R-P) f$ and $2R\sim P+Q$, then $S'$ is
indecomposable, $\e'\sim 0$ and $Y_0=D_R'$.

\item If $x\in D_R\cap D_T$ with $R\neq T$ and $R+T\sim P+Q$, then $S'$ is
decomposable, $\e'\sim 2R-P-Q$ and $Y_0=D_R'$.

\end{enumerate}

\end{prop}
{\bf Proof:}

By Proposition \ref{irreduciblesnodescomp1} we know that two curves $D_Q\sim X_0+(R-P)
f$ and $D_T\sim X_0+(T-P)f$ satisfying $R+T\sim P+Q$ pass through a generic point of
the generator $Qf$; if $T=R$, then it passes a unique curve. Hence, we have two cases: 

\begin{enumerate}

\item If it passes a unique curve $D_R$ through $x$, then $D_R'^2=D_R^2-1=0$. Any
other curve $D_T$ of minimum self-intersection does not pass through $x$ so
$D_T'^2=D_T^2+1=2$. Other curve $D$ of $S$ is in a linear system $|X_0+\b f|$ with
$deg(\b)\geq 1$, so $D^2\geq 3$ and $D'^2\geq D^2-1\geq 1$. From this we deduce that
the minimum self-intersect curve of $S'$ is $D_R'$ and this meets any other curve:
$D_R'.D'=(D_R'^2+D'^2)/2\geq 1$. It follows that $S'$ is indecomposable and $\e'\sim
D_R'^2\sim (X_0+(R-P)f)^2+Q\sim 2R-P+Q\sim 0$.

\item If $x\in D_R \cap D_T$, then, applying elementary transformation
properties, we have $D_R'^2=D_T'^2=0$ and $D_R'.D_T'=D_R.D_T-1=0$. Hence, the ruled
surface $S'$ is decomposable, $D_R'$ is one of the minimum self-intersection curves of
$S'$ and $\e'\sim D_R'^2\sim 2R-P+Q\not\sim 0$. \qed

\end{enumerate}

\bigskip

In the above proposition we have seen that any indecomposable elliptic ruled surface
with invariant $e=0$ or $e=-1$ is obtained by applying an elementary transformation to a
decomposable one.

Moreover, if we apply an elementary transformation to any of two models, we do not
obtain new models of indecomposable surfaces.

Thus, we have actually proved the following theorem:

\begin{teo}\label{elipticasnodescomp}

There only exist two models of indecomposable elliptic ruled surfaces. They have
invariants $e=0$ and $e=-1$.

\end{teo}

\begin{rem}\label{notaelipticanodescomp1}

{\em The indecomposable ruled surface $S$ with invariant $e=-1$ and $\e\sim P$ is obtained
by projecting a decomposable one from a generic point on a generator. If we change the
generator, $\e$ is modified too.

However, both models are isomorphic. We have a family of minimum self-intersection
curves on $S$ parameterized by $X$. If $\e\sim P$, then $\Te_{X_0}(X_0)\sim \Te_{X}(P)$.
Taking other curve $D_R\sim X_0+(R-P)f$, we have $\Te_{D_R}(D_R)\cong \Te_{X}(2R-P)$.
If we take $R$ satisfying $2R\sim P+Q$, then $\Te_{D_R}(D_R)\cong \Te_{X}(Q)$. So, by
considering $D_Q$ as a minimum self-intersection curve, we obtain $\e\sim Q$.

Note that changing minimum self-intersection curve corresponds to modifying the
normalized model of $S$. If we have the normalized sheaf $\E_0$ with $\e\sim P$, when we
consider $D_R$ as minimum self-intersection curve, we are taking a new normalized
sheaf $\E_0'=\E_0\otimes \Te_X(R-P)$.}

\end{rem}

\bigskip

\begin{rem}\label{notanagata}
Nagata Theorem asserts that every geometrically ruled surface $\pi:\P(\E_0)\lrw X$
 is obtained from $X\times \P^1$ by applying a finite number of elementary 
transformations.

We can study this fact in the elliptic ruled surfaces. Let us remember that $X\times
P^1$ corresponds to the ruled surface $\P(\Te_X\oplus \Te_X)$. By applying
Theorem~\ref{telementalelipticasdescomp} we can indicate with detail which is the minimum number of
elementary transformations that we need to obtain each model of elliptic ruled surface:

\begin{enumerate}

\item The decomposable elliptic ruled surface with $e=0$ and $\e\not\sim 0$ is
obtained by applying $2$ elementary transformations to $X\times P^1$ in $2$ generic points.

\item The decomposable elliptic ruled surface with $e=1$ is
obtained by applying $1$ elementary transformation to $X\times P^1$ in a generic point.

\item The decomposable elliptic ruled surface with $e>1$  is
obtained by applying $e$ elementary transformations to $X\times P^1$ in $e$ points laying in
the same directrix curve $X\times \{x_0\}$.

\item The indecomposable elliptic ruled surface with $e=0$ and $\e\not\sim 0$ is
obtained by applying $2$ elementary transformations to $X\times P^1$ in $2$ infinitely near
points: the second elementary transformation is applied in a generic point of the
exceptional divisor corresponding to the first one.

\item The indecomposable elliptic ruled surface with $e=-1$ is
obtained by applying $3$ elementary transformations to $X\times P^1$ in $3$ generic points.

\end{enumerate}

\end{rem}

\bigskip

Let us now study base-point-free linear systems on indecomposable elliptic ruled
surfaces. In this way, we will describe indecomposable elliptic
scrolls.

\begin{prop}\label{pfijoselipticasnodescomp}

Let $S$ be an indecomposable elliptic ruled surface and let $|H|=|X_0+\b f|$ be a
complete linear system on $S$. $|H|$ is base-point-free if and only if $deg(\b)\geq
2+e$.

\end{prop}
{\bf Proof:}

There are two cases, with $e=0$ or $e=-1$. In both of them, if $deg(\b)\geq 2$, then
$\b$ is nonspecial and $\b$ and $\b+\e$ are base-point-free. By applying Corollary $1.8$ in
(\cite{fuentes}, we deduce that the linear system $|H|$ is base-point-free.
If $deg(\b)\leq 1$ we treat each case independently:

\begin{enumerate}

\item Let us suppose $e=0$.

If $deg(\b)\leq 0$, then, since $h^0(\Te_S(X_0+\b f))\leq
h^0(\Te_X(\b))+h^0(\Te_X(\b+\e))$ and $h^0(\Te_S(X_0))=1$, we have $h^0(\Te_S(H))\leq 1$ and the linear
system has base points.

If $deg(\b)=1$, then $\b+\e\sim \b$ has a base point and, according to Proposition $1.5$ in
(\cite{fuentes}, we see that linear system $|H|$ has a base point.

\item Let us suppose $e=-1$.

If $deg(\b)\leq 0$, then we saw that the linear system $|X_0+\b f|$ has at most a unique
curve so it has base points.

If $deg(\b)=1$, then we know that $\b$ is nonspecial, so $h^0(\Te_S(X_0+P f))=
h^0(\Te_X(P))+h^0(\Te_X(P+\e))=3$. Moreover, we know $h^0(\Te_S(X_0+(\P-Q)))=1$ for
all $Q$. Thus, by Corollary $1.4$ \cite{fuentes}, $|H|$ is base-point-free. \qed

\end{enumerate}

\begin{prop}\label{amplitudelipticasnodescomp}

Let $S$ be an indecomposable elliptic ruled surface and let $|H|=|X_0+\b f|$ a
complete linear system on $S$. $|H|$ is very ample if and only if $deg(\b)\geq 3+e$.

\end{prop}
{\bf Proof:}

According to Theorem $1.10$ in \cite{fuentes}, we know that the linear system $|H|$ is
very ample if and only if $|H|$ and $|H-Pf|$ are base-point-free for all $P\in X$. By
the above theorem, this happens when $deg(\b)\geq e+2$ and $deg(\b-P)\geq e+2$; that is,
when
$deg(\b)\geq e+3$. \qed

\begin{lemma}\label{normalidadelipticasnodescomp}

Let $S$ be an indecomposable elliptic ruled surface and let $|H|=|X_0+\b f|$ be a
complete base-point-free linear system on $S$. $|H|$ defines a regular map
$\phi:S\lrw \P^N$. Let $D$ be an unisecant curve on $S$, with $D\sim X_0+\aa f$. The
curve $\phi(D)$ is linearly normal if and only if $deg(\aa)\leq deg(\b)$ and $\aa
\not\sim
\b$.

\end{lemma}
{\bf Proof:}

Let us consider the trace of $|H|$ on the curve $D$:
$$\begin{array}{l}
0\rw H^0(\Te_S((\b -\aa) f))\lrw
H^0(\Te_S(X_0+\b f))\stackrel{\A}{\lrw} H^0(\Te_{D}(X_0+\b
f))=\\
=H^0(\Te_X(\aa+\b+\e))
\end{array}
$$
The curve $\phi(D)$ is linearly normal when $|H|$ traces on $D$ the complete linear
system $|\aa+\b+\e|$; that is, when $\A$ is a surjection. This happens when
$h^0(\Te_S(X_0+\b f))=h^0(\Te_X(\aa+\b+\e))+h^0(\Te_X(\b-\aa))$. 

Since $|H|$ is base-point-free, by Theorem \ref{pfijoselipticasnodescomp},
$deg(\b)\geq 2+e\geq 1$. Then $\b$ is nonspecial and we have 
$h^0(\Te_S(X_0+\b f))=h^0(\Te_X(\b))+h^0(\Te_X(\b+\e))$.

Moreover, according to Propositions \ref{irreduciblesnodescomp0} and \ref{irreduciblesnodescomp1}, we
know that
$deg(\aa)\geq 0$. From this,
$deg(\aa+\e+\b)\geq 2$ and divisor $\aa+\e+\b$ is nonspecial. By Riemann-Roch,
$h^0(\Te_X(\aa+\b+\e))=deg(\aa)+deg(\b)+deg(\e)$ and
$h^0(\Te_X(\b-\aa))=deg(\b)-deg(\aa)+h^1(\Te_X(\b-\aa))$.

As $deg(\b)\geq 2+e$ then $\b$ and $\b+\e$ are nonspecial on $X$ and $h^0(\Te_S(X_0+\b
f))=2deg(\b)-e$. Thus
$$h^0(\Te_S(X_0+\b f))-(h^0(\Te_X(\aa+\b+\e))+h^0(\Te_X(\b-\aa)))=-h^1(\Te_X(\b-\aa)),$$
and $\A$ is a surjection when $\b-\aa$ is nonspecial; equivalently,
when $\b\not\sim \aa$ and $deg(\aa)\leq deg(\b)$. \qed

\begin{teo}\label{scrollselipticasnodescomp}

Let $S$ be an indecomposable elliptic ruled surface. Let $|H|=|X_0+\b f|$ be a
complete base-point-free linear system on $S$. Denote $b=deg(\b)$. Let $\phi:S\lrw
\P^N$ be the regular map defined by $|H|$ on $S$ and let $R=\phi(S)$ be the image scroll.
Then $R$ is one of the following models:

\begin{enumerate}

\item When $\phi$ is not birational:

\begin{enumerate}

\item If $e=-1$ and $b=1$, then $R$ is a plane and $\phi$ is a $3:1$ map from $S$
onto $R$.

\end{enumerate}

\item When $\phi$ is birational:

\begin{enumerate}

\item If $e=0$ and $b=2$, then $R$ is a elliptic scroll of degree $4$ in  $\P^3$. It is
generated by a $1:2$ correspondence with an united point between a line and a
linearly normal elliptic curve meeting at a point.

The singular locus of $R$ is the double line.

$R$ has unisecant elliptic curves of degree $d\geq 3$ in the linear systems $|X_0+\aa f|$
with $deg(\aa)=d-2$. These curves are linearly normal if and only if $d\leq 4$ and
$\aa\not\sim \b$ (if $\aa\sim \b$ they are the hyperplane sections).

\item If $e=-1$ and $b\geq 2$, then $R$ is a smooth elliptic scroll of degree $2b+1$
in $\P^{2b}$. It is generated by a $1:1$ correspondence with an united point between
two linearly normal elliptic curves of degree $b+1$ meeting at a point.

$R$ has unisecant elliptic curves of degree $d\geq b+1$ in the linear systems $|X_0+\aa f|$
with $deg(\aa)=d-b+1$. These curves are linearly normal if and only if $d\leq 2b+1$ and
$\aa\not\sim \b$ (if $\aa\sim \b$ they are the hyperplane sections).

\item If $e=0$ and $b\geq 3$, then $R$ is a smooth elliptic scroll of degree $2b$ in 
$\P^{2b-1}$. It is generated by a $1:1$ correspondence with an united point between two
linearly normal elliptic curves of degrees $b$ and $b+1$ meeting at a point.

There is a unique unisceant curve of minimum degree $b$. Moreover, $R$ has unisecant
elliptic curves of degree $d\geq b+1$ in the linear systems
$|X_0+\aa f|$ with $deg(\aa)=d-1$. These curves are linearly normal if and only if $d\leq
2b$ and $\aa\not\sim \b$ (if $\aa\sim \b$ they are the hyperplane sections).

\end{enumerate}

\end{enumerate}

\end{teo}
{\bf Proof:}

We begin by studying what happens when $|H|$ is base-point-free but not very ample. By
Proposition \ref{pfijoselipticasdescomp}, $|H|$ is base-point-free when $deg(\b)\geq
e+2$.

Propositions \ref{irreduciblesnodescomp0} and \ref{irreduciblesnodescomp1} determine us the
families of irreducible curves and by Lemma \ref{normalidadelipticasnodescomp} we
know when they are linearly normal.

Note that the degree of $R$ is $(X_0+\b f)^2=2b-e$. $R$ lies in $\P^N$ with 
$N=h^0(\Te_S(X_0+\b f))$. Since $b\geq e+2\geq 1$, $\b$ is nonspecial,
$h^0(\Te_S(X_0+\b f))=h^0(\Te_X(\b))+h^0(\Te_X(\b+\e))$ and $N=2b-e-1$.

In order to find the singular locus of the scroll, we will use Theorem $1.10$ in
\cite{fuentes}, $\phi$ is not an isomorphism at the base points of the linear systems $|H-P
f|$.

\begin{enumerate}

\item If $e=-1$ and $b=1$ then $|H|=|X_0+Pf|$. For any point $Q\in X$, $|X_0+(P-Q)f|$
has a unique curve $D_Q$. Their points are base points of $|X_0+(P-Q)f|$. Moreover,
by Proposition \ref{irreduciblesnodescomp1}, the curves $D_Q$ fill the surface, so all
points of $S$ are base points for some system $|H-Qf|$ and $\phi$ is not birational.
In fact, since $h^0(\Te_S(X_0+P f))=3$, $R$ is a plane. Finally, $(X_0+Pf)^2=3$ so
the morphism $\phi$ is a $3:1$ map from $S$ onto a plane.
 
\item Let us suppose $e=0$ and $b=2$. $R$ is a scroll of degree $4$ in $\P^3$. In
order to find the singular locus of $R$, we study the base points of the systems $|H-Pf|$.

Since $b=2$, $\b-P$ is nonspecial, so
$$h^0(\Te_S(H-Pf))=h^0(\Te_X(\b-P))+h^0(\Te_X(\b+\e-P))=2.$$
Given $Q\in X$ with
$P+Q\not\sim \b$, $\b-P-Q$ is nonspecial and $h^0(\Te_S(H-(P+Q)f))=0$. By applying
Proposition $1.5$ in \cite{fuentes}, we see that $|H-Pf|$ have not base points on $Qf$.

If $Q\in X$, but $P+Q\sim \b$, then $|H-Pf|=|X_0+Qf|$. By Proposition~$1.5$ in \cite{fuentes},
 the linear system $|H-Pf|$ has a base point at $X_0\cap Qf$. Since
$h^0(\Te_S(X_0+Q)f))=1$ this base point is unique on $Qf$. Moving $P$ on $X$, we see
that
$\b-P\sim Q$ becomes any point of $X$. So all points of $X_0$ are base points for some
system $|X_0+(\b -P)f|$.

It follows that $\phi$ is not an isomorphism at points of $X_0$.  The singular locus
of $R$ is $\phi(X_0)$, which is given by the complete linear system $|\b|$ on $X$.
Since $b=2$, $\phi(X_0)$ is a double line.

Finally, by applying Proposition \ref{irreduciblesnodescomp0} and Lemma
\ref{normalidadelipticasnodescomp}, we know that there is a linearly normal
elliptic curve of degree $3$ meeting $X_0$ at one point. Then $R$ is generated by a
$1:2$ correspondence with a united point between line $\phi(X_0)$ and the elliptic
curve.

\item If $e=-1$ and $b\geq 2$, the linear system is very ample. Then, $R$ is a
nonsingular elliptic scroll of degree $2b+1$ in $\P^{2b}$.

By applying Proposition \ref{irreduciblesnodescomp1} and Lemma 
\ref{normalidadelipticasnodescomp}, we see that there are two linearly normal elliptic
curves of degree $b+1$ meeting at one point. The scroll is generated by a $1:1$
correspondence with a united point between these curves. (Actually, we know that there
is an unidimensional family of these curves parameterized by $X$).

\item If $e=0$ and $b\geq 3$, the linear system is very ample. Then, $R$ is a
nonsingular elliptic scroll of degree $2b$ in $\P^{2b-1}$.

By applying Proposition \ref{irreduciblesnodescomp0} and Lemma 
\ref{normalidadelipticasnodescomp}, we see that there are two linearly normal elliptic
curves of degrees $b$ and $b+1$ meeting at one point. The scroll is generated by a
$1:1$ correspondence with a united point between these curves. \qed

\end{enumerate}

\bigskip

Now, we present some tables where all the elliptic scrolls of $\P^N$ are described. We explain with
detail their projective generation and singular loci. We indicate the degree of the divisor $\b$
providing the linear system of hyperplane sections $|H|=|X_0+\b f|$. 

\begin{center}
{\footnotesize
\begin{tabular}{|c|c|l|l|c|}
\hline
\hline\multicolumn{5}{|c|}{TABLE 1. ELLIPTIC SCROLLS IN $\P^3$.} \\
\hline
\hline
{$e=-\gr(\e)$}&{$\gr(\b)$}&{Irreducible elements.}&{Projective
generation.}&{Sing.}\\
\hline
\hline
{$e=0$}&{$2$}&{$|X_0|$}&{$C^2_1\stackrel{(1:2)}{\lrw}C^3_1$}&{$C^2_1$}\\
{$\e\sim 0$}&{}&{$|X_0+\aa f|$}&{$1$ united point}&{}\\
{}&{}&\multicolumn{1}{|r|}{$deg(\aa)\geq 1$}&{Degree $4$.}&{}\\
\hline
{$\e\sim 0$}&{$2$}&\multicolumn{3}{|c|}{ Degenerated case. Double quadric.}\\
\hline
{$e=0$}&{$2$}&{$|X_0|$}&{$C^2_1\stackrel{(2:2)}{\lrw}C^2_1$}&{$C^2_1,C^2_1$}\\
{$\e\sim P-Q$}&{}&{$|X_1|$}&{Degree $4$.}&{}\\
{}&{}&{$|X_0+\aa f|$}&{}&{}\\
{}&{}&\multicolumn{1}{|r|}{$deg(\aa)\geq 1$}&{}&{}\\
\hline
{$e=3$}&{$3$}&{$|X_0|$}&{Cone over $C^3_1$}&{$V$}\\
{}&{$\b\sim -\e$}&{$|X_1|$}&{and vertex $V(=X_0)$.}&{}\\
{}&{}&{$|X_0+\aa f|$}&{Degree $3$.}&{}\\
{}&{}&\multicolumn{1}{|r|}{$deg(\aa)\geq 4$}&{Speciality $1$.}&{}\\
\hline
\end{tabular}}
{\footnotesize
\begin{tabular}{|c|c|l|l|c|}
\hline
\hline\multicolumn{5}{|c|}{TABLE 2. ELLIPTIC SCROLLS IN $\P^N$ ($N$ odd, $N\geq 5$). } \\
\hline
\hline
{$e=-\gr(\e)$}&{$\gr(\b)$}&{Irreducible elements.}&{Projective
generation.}&{Sing.}\\
\hline
\hline
{$e=0$}&{$\frac{N+1}{2}$}&{$|X_0|$}&
{$C^{\frac{N+1}{2}}_1\stackrel{(1:1)}{\lrw}C^{\frac{N+1}{2}+1}_1$}&{$\emptyset$}\\
{$\e\sim 0$}&{}&{$|X_0+\aa f|$}&{$1$ united point}&{}\\
{}&{}&\multicolumn{1}{|r|}{$deg(\aa)\geq 1$}&{Degree $N+1$.}&{}\\
\hline
{$e=0$}&{$\frac{N+1}{2}$}&{$|X_0|$}&
{$C^{\frac{N+1}{2}}_1\stackrel{(1:1)}{\lrw}C^{\frac{N+1}{2}}_1$}&{$\emptyset$}\\
{$\e\sim 0$}&{}&{$|X_0+\aa f|$}&{Degree $N+1$.}&{}\\
{}&{}&\multicolumn{1}{|r|}{$deg(\aa)\geq 2$}&{}&{}\\
\hline
{$e=0$}&{$\frac{N+1}{2}$}&{$|X_0|$}&
{$C^{\frac{N+1}{2}}_1\stackrel{(1:1)}{\lrw}C^{\frac{N+1}{2}}_1$}&{$\emptyset$}\\
{$\e\sim P-Q$}&{}&{$|X_1|$}&{Degree $N+1$.}&{}\\
{}&{}&{$|X_0+\aa f|$}&{}&{}\\
{}&{}&\multicolumn{1}{|r|}{$deg(\aa)\geq 1$}&{}&{}\\
\hline
{$2\leq e<N-3$}&{$\frac{N+1+e}{2}$}&{$|X_0|$}&
{$C^{\frac{N+1-e}{2}}_1\stackrel{(1:1)}{\lrw}C^{\frac{N+1+e}{2}}_1$}&{$\emptyset$}\\
{$e$ even}&{}&{$|X_1|$}&{Degree $N+1$.}&{}\\
{}&{}&{$|X_0+\aa f|$}&{}&{}\\
{}&{}&\multicolumn{1}{|r|}{$deg(\aa)\geq e+1$}&{}&{}\\
\hline
{$e=N-3$}&{$\frac{N+1+e}{2}$}&{$|X_0|$}&
{$C^2_1\stackrel{(1:2)}{\lrw}C^{N-1}_1$}&{$C^2_1$}\\
{$e$ even}&{}&{$|X_1|$}&{Degree $N+1$.}&{}\\
{}&{}&{$|X_0+\aa f|$}&{}&{}\\
{}&{}&\multicolumn{1}{|r|}{$deg(\aa)\geq N-2$}&{}&{}\\
\hline
{$e=N$}&{$N$}&{$|X_0|$}&{Cone over $C^N_1$}&{$V$}\\
{}&{$\b\sim -\e$}&{$|X_1|$}&{and vertex $V(=X_0)$.}&{}\\
{}&{}&{$|X_0+\aa f|$}&{Degree $N$.}&{}\\
{}&{}&\multicolumn{1}{|r|}{$deg(\aa)\geq N+1$}&{Speciality $1$.}&{}\\
\hline
\end{tabular}}

{\footnotesize
\begin{tabular}{|c|c|l|l|c|}
\hline
\hline\multicolumn{5}{|c|}{TABLE 3. ELLIPTIC SCROLLS IN $\P^N$ ($N$ even, $N\geq 4$). } \\
\hline
\hline
{$e=-\gr(\e)$}&{$\gr(\b)$}&{Irreducible elements.}&{Projective
generation.}&{Sing.}\\
\hline
\hline
{$e=-1$}&{$\frac{N}{2}$}&{$|X_0+\aa|$}&
{$C^{\frac{N}{2}+1}_1\stackrel{(1:1)}{\lrw}C^{\frac{N}{2}+1}_1$}&{$\emptyset$}\\
{$\e\sim P$}&{}&\multicolumn{1}{|r|}{$deg(\aa)\geq 0$}&{$1$ united point}&{}\\
{}&{}&{}&{Degree $N+1$.}&{}\\
\hline
{$1\leq e<N-3$}&{$\frac{N+1+e}{2}$}&{$|X_0|$}&
{$C^{\frac{N+1-e}{2}}_1\stackrel{(1:1)}{\lrw}C^{\frac{N+1+e}{2}}_1$}&{$\emptyset$}\\
{$e$ odd}&{}&{$|X_1|$}&{Degree $N+1$.}&{}\\
{}&{}&{$|X_0+\aa f|$}&{}&{}\\
{}&{}&\multicolumn{1}{|r|}{$deg(\aa)\geq e+1$}&{}&{}\\
\hline
{$e=N-3$}&{$\frac{N+1+e}{2}$}&{$|X_0|$}&
{$C^2_1\stackrel{(1:2)}{\lrw}C^{N-1}_1$}&{$C^2_1$}\\
{$e$ odd}&{}&{$|X_1|$}&{Degree $N+1$.}&{}\\
{}&{}&{$|X_0+\aa f|$}&{}&{}\\
{}&{}&\multicolumn{1}{|r|}{$deg(\aa)\geq N-2$}&{}&{}\\
\hline
{$e=N$}&{$N$}&{$|X_0|$}&{Cone over $C^N_1$}&{$V$}\\
{}&{$\b\sim -\e$}&{$|X_1|$}&{and vertex $V(=X_0)$.}&{}\\
{}&{}&{$|X_0+\aa f|$}&{Degree $N$.}&{}\\
{}&{}&\multicolumn{1}{|r|}{$deg(\aa)\geq N+1$}&{Speciality $1$.}&{}\\
\hline
\end{tabular}}

\end{center}

\bigskip

\section{$2$-secant linear systems on an elliptic ruled surface.}\label{dossecantes}

We study the families of $2$-secant curves on an elliptic ruled surface. In order to
get this, we work with the linear systems $|2X_0+\b f|$. We investigate when are they base-point-free
and very ample. Then, we apply Bertini Theorems to determine when is the generic element 
irreducible. 
We begin by working with the decomposable elliptic ruled surfaces.

\begin{prop}\label{pfijoselipticasdescomp2}

Let $S$ be a decomposable elliptic ruled surface and let $|H|=|2X_0+\b f|$ be a complete $2$-secant
linear system. Then, $|H|$ is base-point-free if and only if:

\begin{enumerate}

\item $\b\sim -2\e$ or $deg(\b)\geq 2e+2$, when $e>0$.

\item $\b\sim 0$ or $deg(\b)\geq 2$, when $e=0$ and $\e\sim 0$. 

\item $\b\sim 0$ and $2\e\sim 0$ or $deg(\b)\geq 2$, when $e=0$ and $\e\not\sim 0$. 

\end{enumerate}

\end{prop}
{\bf Proof:}
We use Proposition $2.11$ in \cite{fuentes}. The linear system $|H|$ is
base-point-free when $\b$ and
$\b+2\e$ are base-point-free:

\begin{enumerate}

\item If $e>0$, then $\b+2\e$ is base-point-free when $\b\sim -2\e$ or $\b\geq 2e+2$. In both cases,
$\b$ is base-point-free too, because $e>0$.

\item Let us suppose $e=0$ and $\e\sim 0$. Then, it is sufficient that $\b$ is base-point-free, that
is, $\b\sim 0$ or $\b\geq 2$.

\item Let us suppose $e=0$ and $\e\not\sim 0$. Then, if $2\e\sim 0$ we are in the above situation.
When $2\e\not\sim 0$, if $\b\sim 0$, then $\b+2\e\not\sim 0$ and it has base points; and conversely.
Thus $\b$ and $\b+2\e$ are base-point-free when $b\geq 2$. \qed

\end{enumerate} 

\begin{prop}\label{amplitudelipticasdescomp2}

Let $S$ be a decomposable elliptic ruled surface and let $|H|=|2X_0+\b f|$ be a complete $2$-secant
linear system. Then, $|H|$ is very ample if and only if $deg(\b)\geq 2e+3$.

\end{prop}
{\bf Proof:}

By Theorem $2.13$ in \cite{fuentes}, we know that the linear system is very ample
if and only if $\b$ and $\b+2\e$ are very ample, that is, if $deg(\b)\geq 3$ and $deg(\b)\geq 2e+3$.
Since $e\geq 0$, it is sufficient that $\b\geq 2e+3$. \qed

\begin{prop}\label{irreduciblesdescomp2}

Let $S$ be a decomposable elliptic ruled surface and let $|H|=|2X_0+\b f|$ be a complete $2$-secant
linear system. Then the generic element of $|H|$ is irreducible if only if it satisfies one of the
following conditions:

\begin{enumerate}

\item $deg(\b)\geq 2e+2$.

\item $deg(\b)=2e+1$ and $\e\not\sim 0$.

\item $\b\sim -2\e$ and $e>0$.

$\b\sim -2\e$, $e=0$, $\e\not\sim 0$ and $2\e\sim 0$.

\end{enumerate}

Moreover, if the generic element is irreducible then it is smooth too and its genus is $deg(\b+\e)+1$.

\end{prop}
{\bf Proof:}

We can estimate the genus by using the formula 
$$g(C+D)=g(C)+g(D)+CD-1$$
by using (\cite{hartshorne}, V, ex. 1.3.) and by taking $C=X_0$ and $D=X_0+\b f$.

If $D\sim 2X_0+\b f$ is irreducible, the divisors $\pi_*(D\cap X_0)\sim \b+2\e$ and $\pi_*(D\cap
X_1)\sim \b$ must be effective, so necessarily $deg(\b)\geq 2e$.

Let $\vhi_{|H|}:S\lrw \P^N$ be the rational map defined by the linear system $|H|$:

\begin{enumerate}

\item If $deg(\b)\geq 2e+3$, the linear system $|H|$ is very ample. So the generic element is
irreducible and smooth.

\item If $deg(\b)=2e+2$, the linear system $|H|$ is base-point-free. In particular, $\b$,
$\b+\e$ and
$\b+2\e$ are base-point-free. By Proposition $2.12$ in \cite{fuentes}, the map
$\vhi_{|H|}$ applies the generators onto nonsingular conics. Moreover, the image of $X_0$ is given by
the complete linear system $|\b+2\e|$. Since it has degree $2$, it is a double line. Then
$dim(Im(\vhi_{|H|}))=2$ and by Bertini Theorem, we deduce that the generic element of the linear
system is irreducible and smooth.

\item If $deg(\b)=2e+1$ and $e>0$, by Proposition $2.11$ in \cite{fuentes}, we see
that the linear system
$|H|$ has a unique base point at $X_0\cap Pf$, where $\b+2\e\sim P$. A generic generator $Qf$ is
applied onto a nonsingular conic. The image of $X_1$ is given by the complete linear system
 $|\b|$ on $X$. Since $deg(\b)\geq 3$, $\vhi_{|H|}(X_1)$ is a nonsingular elliptic curve. Thus
$dim(Im(\vhi_{|H|}))=2$ and by Bertini Theorem, the generic element of $|H|$ is irreducible and it
has at most a singular point at  $X_0\cap Pf$. But $HX_0=1$, so the generic element of
 $H$ meets $X_0$ at a unique point and it is smooth at the points of $X_0$.

\item If $deg(\b)=2e+1$, $e=0$ and $\e\not\sim 0$, the linear system has base points at $X_0\cap
P_0f$ and $X_1\cap P_1f$, where $\b+2\e\sim P_0$ and $\b\sim P_1$. Anyway, the generic generator is
applied onto a nonsingular conic. Moreover, $\b+\e\sim P$ and $P\neq P_i$ because
$\e\not\sim 0$. Therefore, $Pf$ is applied onto a double line (see the proof of
(\cite{fuentes}, 2.11)). Thus $dim(Im(\vhi_{|H|}))=2$ and by Bertini Theorem the generic
element of $|H|$ is irreducible and it has at most singular points at $X_i\cap
P_i f$. Since $HX_i=1$, the generic element is smooth at the points of $X_0$ and $X_1$.

\item If $deg(\b)=2e+1$ and $\e\sim 0$, then $h^0(\Te_S(2X_0+\b f))=3=h^0(\Te_S(2X_0))$, so
$\b f$ is a fixed component of the linear system and the generic element is reducible.

\item If $deg(\b)=2e$, then necessarily $\b\sim -2\e$, because $\b+2\e$ must be an effective divisor.

If $e>0$, the linear system is base-point-free. The generators are
applied onto nonsingular conics. The curve $X_1$ is applied onto an elliptic curve of degree $\b\geq
2$, so $dim(Im(\vhi_{|H|}))=2$ and by Bertini Theorem, the generic element of the system is
irreducible and smooth.

If $e=0$ and $e\not\sim 0$, since $\b$ must be and effective divisor, $-2\e\sim 0$. In this case the
linear system is base-point-free. By Bertini Theorem the generic element is smooth. A smooth
element in $|H|$ does not contain generators. If it is reducible, it has two disjoint unisecant
curves. They must be $X_0+X_1$. Since $h^0(\Te_S(2X_0-2\e f))=2$, the reducible elements don't fill
the linear system. Thus the generic element is irreducible and smooth.

If $e=0$ and $\e\sim 0$, then $h^0(\Te_S(2X_0))=3$ . Since $h^0(\Te_S(X_0))=2$, $|H|=\{D+D'/D,D'\sim
X_0\}$, so the generic element is reducible. In fact, $Im(\vhi_{|H|})$ is a conic whose
hyperplane sections parameterize the curves of $|2X_0|$. \qed

\end{enumerate}

Now, let us study the $2$-secant linear systems in the indecomposable ruled surfaces. First, we will
generalize some results about $m$-secant divisors on decomposable ruled surfaces that appear in
(\cite{fuentes}, 2).

\begin{lemma}\label{telementaldimsistm}
Let $\pi:S\lrw X$ be a geometrically ruled surface and let $\nu:S'\lrw S$ be the elementary
transformation at the point $x\in S$, $x\in Pf$. Let $C$ be a
$m$-secant and $\aa$ a divisor on $X$. Then:

\begin{enumerate}

\item $|\nu^*(C)+\aa f)|\cong |C+(\aa+mP) f-mx|$.

\item $h^0(\Te_{S'}(\nu^*(C)+\aa f))=h^0(\Te_{S}(C+(\aa+mP) f-mx))$. \qed

\end{enumerate}

\end{lemma}
{\bf Proof:}

Let $|C+\b f|$ be a $m$-secant complete linear system in $S$. Let $D$ be a curve of the linear
system. Then $\nu^*(D)\sim \nu^*(C)+\b f$ and $D'+\mu_x(D)Pf\sim \nu^*(C)+\b f$. From this, $D'\sim
\nu^*(C)+(\b-\mu_x(D)P)f$ and the elements of the linear system
$|\nu^*(C)+\aa f|$ come from the elements of the linear systems $|C +(\aa+kP)f-kx|$:
$$
h^0(\Te_{S'}(\nu^*(C)+\aa f))=\max\{h^0(\Te_{S}(C+(\aa+kP) f-kx))/k=0,\dots ,m\}
$$
Since $|C+(\aa+ kP)f-kx|\subset |C+(\aa+ (k+1)P)f-(k+1)x|$, the conclusion follows.\qed

\begin{lemma}\label{cotadimsistmsecante}
Let $S$ be a ruled surface and let $|H|=|mX_0+\b f|$ be a $m$-secant linear system on $S$.
Then:
$$
h^0(\Te_S(X_0+\b f))\leq \sum\limits^m_{k=0} h^0(\Te_X(\b + k\e))
$$
Moreover, if $\b,\dots,\b+(m-1)\e$ are nonspecial divisors then the equality holds and
$h^1(\Te_S(mX_0+\b f))=h^1(\Te_X(\b+m\e))$.

\end{lemma}
{\bf Proof:}

The proof is by induction on $m$.

If $m=1$ we consider the exact sequence:
$$
0\lrw \Te_S(\b f)\lrw \Te_S(X_0+\b f)\lrw \Te_{X_0}(X_0+\b f)\lrw 0
$$
By applying cohomology we have:
$$
\begin{array}{rccccccl}
{0}&{\rw}&{H^0(\Te_X(\b f))}&{\rw}&{H^0(\Te_S(X_0+\b f))}&{\rw}&{H^0(\Te_X(\b+\e))}&{\rw}\\
{}&{\rw}&{H^1(\Te_X(\b f))}&{\rw}&{H^1(\Te_S(X_0+\b f))}&{\rw}&{H^1(\Te_X(\b+\e))}&{\rw}\\
{}&{\rw}&{0}&{}&{}&{}&{}&{}\\
\end{array}
$$
Then $h^0(\Te_S(X_0+\b f)\leq h^0(\Te_X(\b))+h^0(\Te_X(\b+\e))$. Moreover, we see that if 
$\b$ is nonspecial the equality holds and $h^1(\Te_S(X_0+\b f))=h^1(\Te_X(\b+\e))$.

Let us suppose that the formula holds for $m-1$. Let $|H|=|mX_0+\b f|$. Consider the exact
sequence:
$$
0\lrw \Te_S(H-X_0)\lrw \Te_S(H)\lrw \Te_{X_0}(H)\lrw 0.
$$
By applying cohomology, we get the long exact sequence
$$
\begin{array}{rccccccl}
{0}&{\rw}&{H^0(\Te_S(H-X_0))}&{\rw}&{H^0(\Te_S(H))}&{\rw}&{H^0(\Te_X(\b+m\e))}&{\rw}\\
{}&{\rw}&{H^1(\Te_S(H-X_0))}&{\rw}&{H^1(\Te_S(H))}&{\rw}&{H^1(\Te_X(\b+m\e))}&{\rw}\\
{}&{\rw}&{0}&{}&{}&{}&{}&{}\\
\end{array}
$$
where $H-X_0\sim (m-1)X_0+\b f$. Then $h^0(\Te_S(mX_0+\b f)\leq h^0(\Te_S((m-1)X_0+\b
f))+h^0(\Te_X(\b+m\e))$. We see that if
$\b,\dots,\b+(m-2)\e$ are nonspecial, then by induction hypothesis
$h^1(\Te_S((m-1)X_0+\b f))=h^1(\Te_X(\b+(m-1)\e))$. Moreover, if $\b+(m-1)\e$ is nonspecial too, then
the equality holds and 
$h^1(\Te_S(mX_0+\b f))=h^1(\Te_X(\b+m\e))$. \qed

\begin{lemma}\label{suficientepfijosmsecante}
Let $S$ be a ruled surface and let $|H|=|mX_0+\b f|$ be a complete $m$-secant linear system. If 
$$
h^0(\Te_S(mX_0+(\b-P) f))=h^0(\Te_S(mX_0+\b f))-(m+1)
$$
then the linear system is base-point-free on the generator $Pf$. Moreover, this is applied on a
linearly normal smooth rational curve of degree $m$ by the rational map $\vhi_{|mX_0+\b f|}:S\lrw
\P^N$.
\end{lemma}
{\bf Proof:}

Let us consider the trace of the linear system $|H|$ on $Pf$:
$$
0\lrw H^0(\Te_S(mX_0+(\b-P)f))\lrw H^0(\Te_S(mX_0+\b f))\stackrel{\A}{\lrw}
H^0(\Te_{\P^1}(m))
$$
If $h^0(\Te_S(mX_0+(\b-P) f))=h^0(\Te_S(mX_0+\b f))-(m+1)$, then the map $\A$ is a surjection and
$|H|$ traces the complete linear system of divisors of degree $1$ in $\P^1$. Thus
$|H|$ is base-point-free on the generator $Pf$ and it is applied onto a linearly normal rational curve
of degree $m$ by the rational map $\vhi_{|mX_0+\b f|}$. \qed

\begin{lemma}\label{pfijomsecante}
Let $S$ be a ruled surface and let $|H|=|mX_0+\b f|$ be a complete $m$-secant linear system:
\begin{enumerate}

\item If $P$ is a base point of $\b+m\e$, then $H$ has a base point at $Pf\cap X_0$.

\item If $|H|$ is very ample then $|\b+m\e|$ is very ample.

\end{enumerate}

\end{lemma}
{\bf Proof:}

Let us consider the trace of the linear system $|H|$ on the curve $X_0$:
$$
H^0(\Te_S(H-X_0))\lrw H^0(\Te_S(H))\lrw H^0(\Te_{X_0}(H))\cong
H^0(\Te_X(\b+m\e))
$$
\begin{enumerate}

\item If $P$ is base point of $\b+m\e$, then all divisors of $|H|$ meet $X_0$
at
$X_0\cap Pf$, so this is a base point of $|H|$.

\item If $|H|$ is very ample, then the rational map $\vhi_{|H|}$ is an isomorphism. In
particular, the restriction $\vhi_{|H|}|_{X_0}$ is an isomorphism, so $\b+m\e$ must be very ample.
\qed

\end{enumerate}

\begin{prop}\label{muyampliosmsecantes}
Let $S$ be a ruled surface and let $|H|=|mX_0+\b f|$ be a complete $m$-secant linear system. Let
$\vhi_{|H|}:S\lrw \P^N$ be the rational map defined by $|H|$. If $\vhi_{|H|}$ is an isomorphism on
the generators and the linear systems $|mX_0+(\b-P)f|$ are base-point-free then the linear system 
$|H|$ is very ample.
\end{prop}
{\bf Proof:}

Since $\vhi_{|H|}$ is an isomorphism on the generators, the linear system $|H|$ is base-point-free. 
Let us see that $|H|$ separates points and tangent vectors.

Let $x,y\in S$, $x\in Pf$, $y\in Qf$:
\begin{enumerate}

\item If $Q=P$, then, since the restriction of the rational map $\vhi_{|H|}$ to
the generators defines an isomorphism, the linear system separates points on the same
generator.

\item If $Q\neq P$, since the linear system $|mX_0+(\b-P) f|$ is base-point-free, there exists a
divisor $D\sim mX_0+(\b-P)f$ such that $y\not\in D$. Thus $D+Pf\sim H$, $x\in D+Pf$, but
$y\not\in D+Pf$.

\end{enumerate}

Let $x\in S$, $x\in Pf$, $t\in T_X(S)$:
\begin{enumerate}

\item If $t\in T_x(Pf)$, then there is a divisor which meets $Pf$ at $x$ transversally,
because the restriction of the rational map $\vhi_{|H|}$ to the generators is an isomorphism. 

\item If $t\not\in T_x(Pf)$, since the linear system $|mX_0+(\b-P) f|$ is base-point-free,
there exists a divisor $D\sim mX_0+(\b-P)f$ such that $x\not\in D$. Then, $D+Pf\sim H$, $x\in D+Pf$,
but $t\not\in T_x(Pf)=T_x(D+Pf)$. \qed

\end{enumerate}

Now, we restrict our attention to the study of $2$-secant linear systems on the indecomposable elliptic
ruled surface with $e=0$. Let us remember that it is obtained by applying a elementary
transformation $\nu^*:S\lrw S_0$ at the point $x\in Pf$, $x\not\in X_0\cup X_1$ to the
surface $S_0=\P(\Te_X\oplus \Te_X(-P))$.

\begin{prop}\label{dimsistnodescomp02}
Let $S$ be the indecomposable elliptic ruled surface with $e=0$. Let
$|H|=|2X_0+\b f|$ be an $m$-secant linear system on $S$.
\begin{enumerate}

\item If $deg(\b)\geq 1$, then $h^0(\Te_S(2X_0+\b f))=3deg(\b)$.

\item If $deg(\b)=0$ and $\b\not\sim 0$ or $deg(\b)\leq 0$, then $h^0(\Te_S(2X_0+\b f))=0$.

\item If $\b\sim 0$, then $h^0(\Te_S(2X_0+\b f))=h^0(\Te_S(2X_0))=1$.

\end{enumerate}
\end{prop}
{\bf Proof:}

\begin{enumerate}

\item If $deg(\b)\geq 1$, then $\b$ and $\b+\e$ are nonspecial. We can apply
Lemma~\ref{cotadimsistmsecante} and we obtain that $h^0(\Te_S(2X_0+\b f))=3h^0(\Te_X(\b))=3deg(\b)$.

\item If $deg(\b)=0$ and $\b\not\sim 0$, $\b$ is nonspecial and by Lemma \ref{cotadimsistmsecante}
$h^0(\Te_S(2X_0+\b f))=0$.

If $deg(\b)<0$, then the inequality $h^0(\Te_S(2X_0+\b f))\leq 3h^0(\Te_X(\b))$ holds, where
$3h^0(\Te_X(\b))=0$.

\item If $\b\sim 0$, since $S$ is the elementary transformation of $S_0$ and by
Lemma~\ref{telementaldimsistm}, $h^0(\Te_S(2X_0))=h^0(\Te_{S_0}(2X_0+2Pf-2x))$. As we see at the
section
 (\cite{fuentes}, 2) the linear system $|2X_0+2Pf|_{S_0}$ defines the linear subsystem of
divisors of degree $2$ of $\P^1$ generated by the homogeneous polynomials $x_0^2$ and $x_1^2$ on $Pf$.
From this, the elements of $|2X_0+2Pf|_{S_0}$ have at most double points in $X_0$ or $X_1$, except
when they contain the generator $Pf$. Thus, since $x\not\in X_0\cup X_1$,
$|2X_0+2Pf-2x|_{S_0}=|2X_0+Pf-x|_{S_0}$. But $x$ is not a base-point of any linear system, because it
does not lie on the curves $X_0$ and
$X_1$. Consequently:
$$
\begin{array}{l}
{h^0(\Te_S(2X_0))=h^0(\Te_{S_0}(2X_0+2Pf-2x))=h^0(\Te_{S_0}(2X_0+Pf-x))}\\
{=h^0(\Te_{S_0}(2X_0+Pf))-1=1}\\
\end{array}
$$ \qed

\end{enumerate}

\begin{prop}\label{pfijosnodescomp02}
Let $S$ be the indecomposable elliptic ruled surface with $e=0$. Let
$|H|=|2X_0+\b f|$ be an $m$-secant linear system on $S$. Then:
\begin{enumerate}

\item The linear system $|H|$ is base-point-free if and only if $deg(\b)\geq 2$.

\item If $deg(\b)=1$, then the linear system $|H|$ has a unique base point at $\b
f\cap X_0$.

\end{enumerate}

\end{prop}
{\bf Proof:}

\begin{enumerate}

\item By Lemma \ref{pfijomsecante}, if the linear system $|H|$ is base-point-free, necessarily
$\b+2\e$ is base-point-free, that is, $deg(\b)\geq 2$ or
$\b\sim 0$. But, if $\b\sim 0$, then $h^0(\Te_X(2X_0))=1$ and the linear system has base points.

Conversely, if $deg(\b)\geq2$, by Lemma
\ref{dimsistnodescomp02} we know that
$h^0(\Te_S(H-Pf))=h^0(\Te_S(H))-3$. Applying Lemma \ref{suficientepfijosmsecante} we deduce that
the linear system is base-point-free and the regular map defined by $|H|$ apply the generators onto
nonsingular conics.

\item Let us suppose $deg(\b)=1$ and let $\b\sim P_0$. Then $h^0(\Te_X(H))=3$ and $h^0(\Te_X(H-Pf))=0$,
except when $P_0=P$. In this case $h^0(\Te_X(H-P_0 f))=1$. Thus, the linear system 
$|H|$ has at most base points in the generator $P_0 f$. Moreover, since $h^0(\Te_X(H-P_0
f))=h^0(\Te_X(H))-2$, the complete linear system
$|H|$ traces a $1$-codimension linear subsystem of the divisors of degree $2$ on $P_0 f$. Therefore
the linear system $|H|$ can have at most a base-point in $P_0 f$. Since
$\b+2\e\sim \b$ has a base point $P_0$ and applying Lemma \ref{pfijomsecante} the conclusion follows.
\qed

\end{enumerate}

\begin{prop}\label{mamplionodescomp02}
Let $S$ be the indecomposable elliptic ruled surface with $e=0$. The $m$-secant linear system
$|H|=|2X_0+\b f|$ is very ample if and only if $deg(\b)\geq 3$.
\end{prop}
{\bf Proof:}

By Lemma \ref{pfijomsecante}, if the linear system is very ample, then $\b+2\e$ is very ample, so
 $deg(\b)\geq 3$.

Conversely, if $deg(\b)\geq 3$, then $h^0(\Te_S(H-Pf))=h^0(\Te_S(H))-3$ for any generator $Pf$ and
$|H-Pf|$ is base-point-free. From Lemmas
\ref{suficientepfijosmsecante} and \ref{muyampliosmsecantes} we deduce that $|H|$ is very ample. \qed

\begin{prop}\label{irreduciblesnodescomp02}
Let $S$ be the indecomposable elliptic ruled surface with $e=0$. The generic element of the $2$-secant
linear system $|H|=|2X_0+\b f|$ is irreducible if and only if
$deg(\b)\geq 1$. Moreover, if the generic element is irreducible then it is smooth and it has genus
 $deg(\b)+1$.
\end{prop}
{\bf Proof:}
We can estimate the genus by using the formula $g(C+D)=g(C)+g(D)+C.D-1$ (\cite{hartshorne}, V, ex.
1.3.) taking $C=X_0$ and $D=X_0+\b f$.

If $deg(\b)\geq 3$ the linear system $|H|$ is very ample, so the generic element is irreducible and
smooth.

If $deg(\b)=2$ the linear system $|H|$ is base-point-free. Since
$h^0(\Te_S(H-Pf))=h^0(\Te_S(H))-3$ for any $P\in X$, $|H|$ defines a regular map which applies
the generators onto smooth conics. Let us see which is the image of $X_0$ by this map:
$$
0\lrw H^0(\Te_S(H-X_0))\lrw H^0(\Te_S(H))\stackrel{\A}{\lrw} H^0(\Te_{X_0}(X_0+\b
f))\cong H^0(\Te_X(\b))
$$
Because $deg(\b)=2$ and $h^0(\Te_S(2X_0+\b f))-h^0(\Te_S(X_0+\b f))=2$, the linear system
$|H|$ traces a complete linear system of degree $2$ on $X$. Thus $X_0$ is applied onto a double
line. $Dim(Im(\vhi_{|H|}))=2$ and by the Bertini Theorem, the generic element is irreducible and
smooth.

If $deg(\b)=1$, then the linear system $|H|$ has a unique base point at $X_0\cap
Pf$, with
$\b\sim P$. The generators (except $Pf$) are applied onto nonsingular conics. The image of
$Pf$ is the projection of a smooth conic from a point, that is, a line. Thus, 
$dim(Im(\vhi_{|H|}))=2$ and by Bertini Theorem, the generic element is irreducible and it has at
most a singular point at 
$X_0\cap Pf$. But $HX_0=1$, so the curves of $|H|$ are smooth at points of $X_0$.

If $deg(b)\leq 0$, then $h^0(\Te_S(H))\leq 1$, so the linear system has base points. \qed

Let us study the $2$-secant linear systems in the indecomposable elliptic ruled surface with $e=-1$.

Let us remember that this surface is obtained by applying a elementary transformation $\nu:S\lrw S_0$
at a point $x\in Qf$, $x\not\in X_0\cup X_1$ to the decomposable elliptic ruled surface
$S_0=\P(\Te_X\oplus \Te_X(\e_0))$ with $deg(\e_0)=0$ and $\e_0\not\sim 0$. Thus $\e\sim \e_0+Q$ and we
will take $\e\sim P_0$.

Moreover, we know that $S$ has an one-dimensional family of curves of mini\-mum self-intersection
parameterized by $X$. Given $Q\in X$ we have the curve $D_Q\sim X_0+(Q-P_0)f$.

\begin{prop}\label{dimsistnodescomp12}
Let $S$ be the indecomposable elliptic ruled surface with $e=-1$ and $\e\sim P_0$. Let $|H|=|2X_0+\b
f|$ a
$2$-secant linear system on $S$.
\begin{enumerate}

\item If $deg(\b)\geq 0$, then $h^0(\Te_S(2X_0+\b f))=3deg(\b)+3$.

\item If $deg(\b)=-1$, we have that:

\begin{enumerate}

\item if $-2\b\not\sim 2P_0$ or $-\b\sim P_0$ then $h^0(\Te_S(2X_0+\b f))=3deg(\b)+3=0$.

\item if $-2\b\sim 2P_0$ and $-\b\neq P_0$ then $h^0(\Te_S(2X_0+\b f))=1$.

\end{enumerate}

\item If $deg(\b)\leq -2$, then $h^0(\Te_S(2X_0+\b f))=0$.

\end{enumerate}
\end{prop}
{\bf Proof:}

\begin{enumerate}

\item If $deg(\b)\geq 0$ and $\b\not\sim 0$, then $\b$ and $\b+\e$ are nonspecial divisors; by
Lemma \ref{cotadimsistmsecante} we obtain that $h^0(\Te_S(2X_0+\b f))=3deg(\b)+3$. 

If $\b\sim 0$, then $|2X_0+\b f|=|2D_Q+2(P_0-Q)|$ for any $Q$ with $2(P_0-Q)\not\sim 0$.
Taking $X_0=D_Q$, we can apply Lemma \ref{cotadimsistmsecante} again and we obtain
$h^0(\Te_S(2X_0+\b f))=3$.

\item Let us suppose $deg(\b)=-1$ and let $\b\sim -P$. By Lemma \ref{telementaldimsistm} we know that
$h^0(\Te_S(2X_0-Pf))=h^0(\Te_{S_0}(2X_0+(2Q-P)f-2x))$ where $S_0$ is the decomposable ruled surface
with $\e_0=P_0-Q$.

If $P\neq P_0$, we take $Q=P$. Then $h^0(\Te_{S_0}(2X_0+Qf))=3$. Moreover, we know that
$|2X_0+Qf|_{S_0}$ traces the linear subsystem of divisors of degree $2$ generated by the polynomials
$\{\lambda_0x_0^2,\lambda_1x_0x_1,\lambda_2x_1^2\}$ on $Qf$,
where $\lambda_i=h^0(\Te_X(Q+i\e_0))-h^0(\Te_X(i\e_0))$. In particular, $\lambda_0=0$ and since
$x\not\in X_0\cup X_1$, there is not curves passing through $x$ with multiplicity $2$ in
$|2X_0+Qf|_{S_0}$, except the reducible elements that contain the generator $Qf$. From this,
$|2X_0+Qf-2x|_{S_0}=|2X_0-x|_{S_0}$. If $2\e_0\not\sim 0$, that is, if $2P_0\not\sim 2Q$, then
$h^0(\Te_{S_0}(2X_0))=1$ and $h^0(\Te_{S_0}(2X_0-x))=0$. If $2\e\sim 0$, that is, if $2P_0\sim 2Q$,
then $h^0(\Te_{S_0}(2X_0))=2$ and because the linear system $|2X_0|_{S_0}$ is base-point-free,
$h^0(\Te_{S_0}(2X_0-x))=1$.

If $P=P_0$, we take $Q\neq P$. Arguing as in the above case, we see that 
$|2X_0+(2Q-P)f-2x|_{S_0}=|2X_0+(Q-P)-x|_{S_0}=|X_0+X_1-x|_{S_0}$; but $h^0(\Te_{S_0}(X_0+X_1))=1$
and since $x\not\in X_0\cap X_1$, we deduce that $h^0(\Te_{S_0}(X_0+X_1-x))=0$.

\item If $deg(\b)\leq 2$, then $|2X_0+\b f|\subset |2X_0-Qf|$ where $Q$ verifies
$2Q\not\sim 2P_0$. By the above discussion, $h^0(\Te_S(2X_0+\b
f))\leq h^0(\Te_S(2X_0-Qf))=0$.
\qed

\end{enumerate}

\begin{prop}\label{pfijosnodescomp12}
Let $S$ be the indecomposable elliptic ruled surface with $e=-1$. Let $|H|=|2X_0+\b f|$ be a
$2$-secant linear system on $S$. The linear system $|H|$ is base-point-free if and only if
$deg(\b)\geq 0$.

\end{prop}
{\bf Proof:}

If $deg(\b)>0$, then $h^0(\Te_S(2X_0+\b f))-h^0(\Te_S(2X_0+(\b-P) f))=3$ for any $P\in
X$; by Lemma \ref{suficientepfijosmsecante} the linear system is base-point-free.

If $deg(\b)=0$, let us suppose that the linear system $|2X_0+\b f|$ has a base point at $x\in Pf$. The
family of curves $\{D_Q\}$ fills the surface, so there exists a curve $D_Q\sim
X_0+(Q-P_0)f$ passing through $x$. Let us consider the trace of the linear system $|H|$ on $D_Q$:
$$
\begin{array}{rcl}
{0}&{\lrw}&{H^0(\Te_S(X_0+(P_0-Q)f))\lrw H^0(\Te_S(2X_0+\b f))\stackrel{\A}{\lrw}}\\
{}&{\stackrel{\A}{\lrw}}&{H^0(\Te_{D_Q}(2X_0+\b f))\cong H^0(\Te_X(\b+2Q))}\\
\end{array}
$$
We have that $dim(Im(\A))=h^0(\Te_S(2X_0+\b f))-h^0(\Te_S(X_0+(\b+P_0-Q) f)=2=h^0(\Te_X(\b+2Q))$.
Since $|\b+2Q|$ is base-point-free, the system can not have base points on $D_Q$, so we get a
contradiction. 

If $deg(\b)<0$, then $h^0(\Te_S(2X_0+\b f))\leq 1$ and the linear system has base points. \qed

\begin{prop}\label{mamplionodescomp12}
Let $S$ be the indecomposable elliptic ruled surface with $e=-1$. The $2$-secant linear system 
$|H|=|2X_0+\b f|$is very ample if and only if $deg(\b)\geq 1$.
\end{prop}
{\bf Proof:}

By Lemma \ref{pfijomsecante} we know that if $|H|$ is very ample then $\b+2\e$ must be very ample,
that is, $deg(\b)\geq 1$. 

Conversely, if $deg(\b)\geq 1$, then $h^0(\Te_S(|H-Pf|))=h^0(\Te_S(|H|))-3$ for any generator
$Pf$; moreover, $|H-Pf|$ is base-point-free. From Lemmas
\ref{suficientepfijosmsecante} and \ref{muyampliosmsecantes} we deduce that $|H|$ is very ample. \qed

\begin{prop}\label{irreduciblesnodescomp12}
Let $S$ be the indecomposable ruled surface with $e=-1$ and $\e\sim P_0$. The generic element of the
$2$-secant linear system
$|H|=|2X_0+\b f|$ is irreducible if and only if $deg(\b)\geq 0$ or
$\b\sim -Q$ with
$2Q\sim P_0$ and $Q\not\sim P_0$ (that is, $-\b$ is one of the three ramification points different
from $P_0$ of the map defined by the divisor $2P_0$ on $X$). Moreover, if the generic element of
$|H|$ is irreducible, then it is smooth and it has genus $2deg(\b)+2$.
\end{prop}
{\bf Proof:} 

The genus follows from the formula $g(C+D)=g(C)+g(D)+CD-1$ (\cite{hartshorne}, V, ex.
1.3.) taking $C=X_0$ and $D=X_0+\b f$.

If $deg(\b)\geq 1$ the linear system $|H|$ is very ample, so the generic element is irreducible and
smooth.

If $deg(\b)=0$ the linear system $|H|$ is base-point-free. It defines a regular map $\vhi$ which
applies the generators onto smooth conics.
Let us see which is the image of the curve $X_0$:
$$
0\lrw H^0(\Te_S(H-X_0))\lrw H^0(\Te_S(H))\stackrel{\A}{\lrw} H^0(\Te_{X_0}(H))\cong
H^0(\Te_X(\b+2P_0))
$$
Since $deg(\b)+2P_0=2$ and $h^0(\Te_S(2X_0+\b f))-h^0(\Te_S(X_0+\b f))=2$, the linear system
$|H|$ traces a complete linear system of degree $2$ in $X_0$. Thus,
$X_0$ is mapped onto a double line. $Dim(Im(\vhi_{|H|}))=2$ and by Bertini Theorem the generic
element of $|H|$ is irreducible and smooth.

If $deg(\b)\leq -1$, then $h^0(\Te_S(2X_0+\b f))=0$ except when $\b\sim -Q$ with $2Q\sim P_0$ and
$Q\not\sim P_0$. In this case $h^0(\Te_S(2X_0+\b f))=1$, that is, the linear system has a unique
curve. If this curve were reducible, it would contain generators or an unisecant curve $D\sim X_0+\aa
f$, with
$deg(\aa)\geq 0$. But $h^0(\Te_S(2X_0+\b f -Pf))=h^0(\Te_S(2X_0+\b f-D))=0$. So the unique curve of
the linear system  $|H|$ is irreducible. Since
$D_QH=1$, $H$ meets each $D_Q$ at a unique point and it can not have singular points.
\qed

\end{document}